\documentstyle[12pt]{amsart}
\input amssym.def
\input amssym.tex
\newtheorem{thm}{Theorem}
\newtheorem{lem}{Lemma}
\newtheorem{prop}{Proposition}
\newtheorem{cor}{Corollary}
\newtheorem{defn}{Definition}
\newtheorem{remark}{Remark}

\newtheorem{ack}{Acknowlegements}
\theoremstyle{definition}
\def\nat{{\Bbb N}}

\def\tee{{\Bbb T}}
\def\A{{\cal A}}

\def\L{{\cal L}}
\def\N{{\cal N}}
\def\P{{\cal P}}
\newcommand{\conv}{\operatorname{conv}}

\renewcommand{\Im}{\operatorname{Im}}

\renewcommand{\Re}{\operatorname{Re}}
\def\chix{\raise.5ex\hbox{$\chi$}}

\newcommand{\LH}{{L^1/H_0^1}}

\begin{document}
\baselineskip=18pt
\title{Some Remarks on the Dunford-Pettis Property}
\author{Narcisse Randrianantoanina}
\address{Department of Mathematics, The University of Texas at Austin, 
Austin, TX 78712-1082}
\email{nrandri@@math.utexas.edu}
\maketitle

\begin{abstract}
Let $A$ be the disk algebra, $\Omega$ be a compact Hausdorff space and $\mu$ 
be a Borel measure on $\Omega$. It is shown that
the dual of $C(\Omega,A)$ has the Dunford-Pettis property.
 This proved in particular that  the spaces $L^1(\mu,
\LH)$ and $C(\Omega,A)$ have the Dunford-Pettis property.
\end{abstract}

\section{Introduction}
Let $E$ be a Banach space, $\Omega$ be a compact Hausdorff space and $\mu$ 
be a finite Borel measure on $\Omega$. 
We denote by $C(\Omega,E)$ the space of all $E$-valued continuous functions 
from $\Omega$ and for $1\le p<\infty$, $L^p(\mu,E)$ stands for the space 
of all (class of) $E$-valued $p$-Bochner integrable functions with its 
usual norm. 
A Banach space $E$ is said to have the Dunford-Pettis property if every 
weakly compact operator with domain $E$ is completely continuous (i.e., takes 
weakly compact sets into norm compact subsets of the range space). 
There are several equivalent definitions.
The basic result proved by Dunford and Pettis in \cite{DP} is that the 
space $L^1(\mu)$ has the Dunford-Pettis property. A.~Grothendieck \cite{GRT} 
initiated the study of Dunford-Pettis property in Banach spaces and showed 
that $C(K)$-spaces have this property. 
The Dunford-Pettis property has a rich history; the survey articles by 
J.~Diestel \cite{D2} and A.~Pe{\l}czy\'nski \cite{PL4} are excellent sources 
of information. 
In \cite{D2}, it was asked if the Dunford-Pettis property can be lifted from 
a Banach $E$ to $C(\Omega,E)$ or $L^1 (\mu,E)$. 
M.~Talagrand \cite{T3} constructed counterexamples for these questions so the 
answer is negative in general. 
There are however some positive results. 
For instance, J.~Bourgain showed (among other things) in \cite{BO2} that 
$C(\Omega,L^1)$ and $L^1(\mu,C(\Omega))$ both have the Dunford-Pettis property; 
K.~Andrews \cite{A1} proved that if $E^*$ has the Schur property then $L^1(\mu,E)$ 
has the Dunford-Pettis property. 
In \cite{SS7}, E.~Saab and P.~Saab observed that if $\A$ is a $C^*$-algebra with 
the Dunford-Pettis property then $C(\Omega,\A)$ has the Dunford-Pettis 
property  and they asked (see \cite{SS7} Question 14, p.389) if a similar 
result holds if one considers the disk algebra $A$. 
In this note we provide a positive  answer to the above question by 
showing that the dual of $C(\Omega,A)$ has the Dunford-Pettis property. 
This implies in particular that both $L^1(\mu,\LH)$ and $C(\Omega,A)$  
have the Dunford-Pettis property. 
Our approach is to study  a ``Random version'' of the minimum norm lifting
from $\LH$ into $L^1$. 

The notation and terminology used and not defined in this note can be 
found in \cite{D1} and \cite{DU}. 

\section{Minimum norm lifting}

Let us beging by fixing some notations. 
Throughout, $m$ denotes the normalized Haar measure on the circle $\tee$. 
The space $H_0^1$ stands for the space of integrable functions on $\tee$ 
such that $\hat f(n) = \int_\tee f(\theta) e^{-in\theta}\,dm(\theta)=0$ 
for $n\le 0$. 

It is a well known fact that $A^* = \LH \oplus_1 M_S(\tee)$ 
where $M_S(\tee)$ is the space of singular measures on $\tee$ (see for 
instance \cite{PL4}). 
Consider the quotient map $q:L^1\to \LH$. 
This map has the following important property: 
for each $x\in \LH$, there exists a unique $f\in L^1$ so that 
$q(f)=x$ and $\|f\|= \|x\|$. 
This fact provides a well-defined map called the minimum norm lifting
$$\sigma :\LH \rightsquigarrow L^1 \text{ such that }
q(\sigma (x)) =x \text{ and } \|\sigma (x)\| = \|x\|\ .$$ 
One of the many important features of $\sigma$ is that it preserves weakly 
compact subsets, namely the following was proved in \cite{PL4}. 

\begin{prop} 
If $K$ is a relatively weakly compact subset of $\LH$ then 
$\sigma (K)$ is relatively weakly compact in $L^1$. 
\end{prop}

Our goal in this section is to extend the minimum norm lifting to certain 
classes of spaces that contains $\LH$. 
In particular we will introduce a random-version of the minimum norm lifting.

First we will extend the minimum norm lifting to $A^*$.

 We define a map
$\gamma: \LH \oplus_1 M_s(\tee) \rightsquigarrow L^1 \oplus_1 M_s(\tee)$ as
follows:
 $$\gamma(\{x,s\})=\{\sigma(x),s\}.$$
Clearly $\gamma$ defines a minimum norm lifting from  $A^*$ into $M(\tee)$.

In order to procede to the next extension,
 we need the following proposition: 

\begin{prop}
Let $\sigma$ and $\gamma$ as above then
\begin{enumerate}
\item[a)]
 $\sigma:\LH\rightsquigarrow L^1$ 
is norm-universally measurable (i.e., the inverse image of every norm 
Borel subset of $L^1$ is norm universally measurable  in 
$\LH$);
 \item[b)] $\gamma: A^* \rightsquigarrow M(\tee)$ is weak*-universally
measurable (i.e, the inverse image of every weak*-Borel subset of $M(\tee)$
is weak*-universally measurable in $A^*$).
\end{enumerate}
\end{prop}

\begin{pf}
For a),
notice that $\LH$ and $L^1$ are Polish spaces (with the norm  
topologies) and so is the product $L^1\times \LH$. 
Consider the following subset of $L^1\times \LH$: 
$$\A = \{ (f,x) ; \, q(f)=x,\, \|f\|= \|x\|\}\ .$$
The set $\A$ is a Borel subset of $L^1\times \LH$. 
In fact, $\A$ is the intersection of the graph of $q$ (which is closed) 
and the subset $\A_1 = \{(f,x) ,\|f\|=\|x\|\}$ which is also closed. 
Let $\pi$ be the restriction  on $\A$  of the second projection of 
$L^1\times \LH$ onto $\LH$. 
The operator $\pi$ is of course continuous and hence $\pi (\A)$ is analytic. 
By Theorem 8.5.3 of \cite{CO}, there exists a universally measurable map 
$\phi :\pi (\A) \to L^1$ whose graph belongs to $\A$. 
The existence and the uniqueness of the minimum norm lifting imply that 
$\pi (\A) = \LH$ and $\phi$ must be $\sigma$.

The proof of b) is done with simmilar argument using the fact $A^*$ and
$M(\tee)$ with the weak* topologies are countable reunion of Polish
spaces and their norms are weak*-Borel measurable. 
The proposition is proved. 
\end{pf} 

Let $(\Omega,\Sigma,\mu)$ be a probability space. 
For a measurable function $f:\Omega\to \LH$, the function $\omega\mapsto 
\sigma (f(\omega))$ ($\Omega\to L^1$) is $\mu$-measurable by Proposition~2. 
We define an extension of $\sigma$ on $L^1(\mu,\LH)$ as follows: 
$$\tilde\sigma :L^1\left( \mu,\LH\right) \rightsquigarrow L^1(\mu,L^1)
\text{ with }
\tilde\sigma (f) (\omega) = \sigma (f(\omega)) \text{ for }\ 
\omega\in\Omega\ .$$ 
The map $\tilde\sigma$ is well defined and $\|\tilde\sigma (f)\|=\|f\|$ 
for each $f\in L^1(\mu,\LH)$. 
Also if we denote by $\tilde q:L^1(\mu,L^1)\to L^1(\mu,\LH)$, 
the map $\tilde q(f)(\omega) = q(f(\omega))$, we get that 
$\tilde q(\tilde\sigma (f)) =f$. 

Similarly if $f: \Omega \to A^*$ is weak*-scalarly measurable, the function
$\omega \mapsto \gamma(f(\omega))$ ($\Omega \to M(\tee)$) is weak*-scalarly
measurable . As above we define $\tilde\gamma$ as follows:

\noindent 
for each measure $G \in M(\Omega,A^*)$, fix $g: \Omega \to A^*$ its
 weak*-density with respect to its variation $|G|$. We define
 $$\tilde\gamma(G)(A)=\text{weak*}-\int_A \gamma(g(\omega))\ d|G|(\omega)
\  \text{for all}\  A \in \Sigma.$$
Clearly $\tilde\gamma(G)$ is a measure and it is easy to check that
 $||\tilde\gamma(G)||=||G||$ (in fact $ |\tilde\gamma(G)|=|G|$). 

The rest of this section is devoted to the proof of
the following result that extends the property of $\sigma$ stated in Proposition~1 to 
$\tilde\sigma$. 

\begin{thm}
Let $K$ be a relatively weakly compact subset of $L^1(\mu ,\LH)$. 
The set $\tilde\sigma (K)$ is relatively weakly compact in $L^1(\mu,L^1)$. 
\end{thm}  

We will need few general facts for the proof. In the sequel, we will identify
(for a given Banach space $F$) the dual of $L^1(\mu,F)$ with the space
$L^\infty(\mu,{F_\sigma}^*)$ of all map $h$ from $\Omega$ to $F^*$ that are 
weak*-scalarly measurable and essentially bounded with the uniform norm (see
\cite{IT}). 

\begin{defn} Let $E$ be a Banach space. A series $\sum\limits_{n=1}^\infty x_n$
in $E$ is said to be weakly unconditionally Cauchy (WUC) if for every
$x^* \in E^*$, the series $\sum\limits_{n=1}^\infty |x^*(x_n)|$ is convergent.
\end{defn}

The following lemma is well known:

\begin{lem}  If $S$ is a relatively weakly compact subset of a Banach space $E$,
then 
for every WUC series $\sum\limits_{n=1}^\infty x_n^*$ in $E^*$,\ 
$\lim\limits_{m \to \infty }x_n^*(x)=0$ uniformly
on $S$.
\end{lem}

The following proposition which was proved in 
\cite{RAN4} is the main ingredient for  the proof of Theorem 1.
 For what follows $(e_n)_n$ denote the unit vector basis of $c_0$ 
and $(\Omega,\Sigma,\mu)$ is a probability  space. 

\begin{prop} 
\cite{RAN4} Let $Z$ be a separable subspace of a real Banach space $E$ and 
$(f_n)_n$ be a sequence of maps from $\Omega$ to $E^*$ that are 
weak*-scalarly  measurable and $\sup_n \|f_n\|_\infty \le 1$. 
Let $a<b$ (real numbers) then: \par 
There exist a sequence $g_n\in \conv\{f_n,f_{n+1},\ldots\}$ measurable 
subsets $C$ and $L$ of $\Omega$ with $\mu(C\cup L)=1$ such that 
\begin{enumerate}
\item[(i)] If $\omega\in C$ and $T\in \L(c_0,Z)$, $\|T\|\le 1$; then for 
each $h_n\in \conv \{g_n,g_{n+1},\ldots\}$, either $\limsup\limits_{n\to\infty} 
\langle h_n(\omega),Te_n\rangle \le b$ or $\liminf\limits_{n\to\infty} 
\langle h_n(\omega),Te_n\rangle \ge a$; 
\item[(ii)] $\omega \in L$, there exists $k\in\nat$ so that for each infinite 
sequence of zeroes and ones $\Gamma$, there exists $T\in\L(c_0,Z)$, 
$\|T\|\le 1$ such that for $n\ge k$, 
\begin{align*}
&\Gamma_n =1 \Longrightarrow \langle g_n(\omega),Te_n\rangle \ge b\\
&\Gamma_n =0 \Longrightarrow \langle g_n(\omega),Te_n\rangle \le a\ .
\end{align*}
\end{enumerate}
\end{prop}

We will also make use of the following fact: 

\begin{lem}
(\cite{PL4}, p.45) 
Let $(U_n)_n$ be a bounded sequence of positive elements of $L^1(\tee)$. 
If $(U_n)_n$ is not uniformly integrable, then there exists a W.U.C. 
series $\sum\limits_{\ell=1}^\infty a_\ell$ in the disc algebra $A$ such that: 
$\limsup\limits_{\ell\to\infty} \sup_n |\langle a_\ell,U_n\rangle|>0$.
\end{lem} 

\begin{pf*}{PROOF OF THEOREM 1}

Assume without loos of generality that $K$ is a bounded subset of $L^\infty (\mu,\LH)$. 
The set $\tilde \sigma (K)$ is a bounded subset of $L^\infty (\mu,L^1(\tee))$. 
Let $|\tilde\sigma (K)| = \{|\tilde\sigma (f)|;f\in K\}$. 
Notice that for each $f\in L^1(\mu,\LH)$, there exists $h\in L^\infty (\mu,
H_\sigma^\infty) = L^1(\mu,\LH)^*$ with $\|h\|=1$ and 
$|\tilde\sigma (f)(\omega)| = \tilde\sigma (f)(\omega). h(\omega)$
(the multiplication of the function $\tilde\sigma(f)(\omega) \in L^1(\tee)$ 
with the function $h(\omega) \in H^\infty(\tee)$) 
for a.e. $\omega\in\Omega$.

Consider $\varphi_n = |\tilde\sigma (f_n)|$ be a sequence of 
$L^1(\mu, L^1(\tee))$ with $(f_n)_n\subset K$ and choose $(h_n)_n\in 
L^\infty (\mu,H_\sigma^\infty)$ so that $\varphi_n(\omega) = \tilde\sigma 
(f_n)(\omega). h_n(\omega)$ $\forall\ n\in \nat$.

\begin{lem}
There exists $\psi_n\in \conv \{\varphi_n,\varphi_{n+1},\ldots\}$ so that 
for a.e. $\omega\in\Omega$, 
$$\lim_{n\to\infty} \langle \psi_n (\omega), Te_n\rangle \text{ exists for 
each } T\in \L(c_0,A)\ .$$
\end{lem} 

To prove the lemma, let $(a(k),b(k))_{k\in\nat}$ be an
 enumeration of all pairs 
of rationals with $a(k)<b(k)$. 
We will apply Proposition~3 successively starting from $(\varphi_n)_n$ for 
$E= C(\tee)$ and $Z=A$. 
Note that Proposition~3 is valid only for real Banach spaces so we will 
separate the real part and the imaginary part. 

Inductively, we construct sequences $(\varphi_n^{(k)})_{n\ge1}$ and 
measurable subsets $C_k$, $L_k$ of $\Omega$ satisfying:
\begin{enumerate}
\item[(i)] $C_{k+1}\subseteq C_k$, \ $L_k\subseteq L_{k+1}$,\  $\mu(C_k\cup L_k)=1$ 
\item[(ii)] $\forall\ \omega\in C_k$ and $T\in \L(c_0,A)$, $\|T\|\le 1$ 
and $j\ge k$, either 
\begin{align*}
&\limsup_{n\to\infty} \Re \langle \varphi_n^{(j)}(\omega), 
Te_n\rangle \le b(k)\ \text{ or}\\ 
&\liminf_{n\to\infty} \Re\langle  \varphi_n^{(j)}(\omega),Te_n
\rangle \ge a(k)
\end{align*}
\item[(iii)] $\forall\ \omega\in L_k$, there exists $\ell\in \nat$ so that 
for each $\Gamma$ infinite sequences of zeroes and ones, there exists 
$T\in \L(c_0,A)$, $\|T\|\le 1$ such that if $n\ge\ell$, 
\begin{align*} 
&\Gamma_n = 1\Rightarrow \Re \langle \varphi_n^{(k)}(\omega),Te_n\rangle
\ge b(k)\\
&\Gamma_n =0 \Rightarrow \Re\langle \varphi_n^{(k)}(\omega),Te_n\rangle 
\le a(k)\ ;
\end{align*}
\item[(iv)] $\varphi_n^{(k+1)}\in \conv \{\varphi_n^{(k)},\varphi_{n+1}^{(k)},
\ldots\}\ .$
\end{enumerate}
Again this is just an application of Proposition 3 starting from the 
sequence $\Omega\to C(\tee)^*$ $(\omega\mapsto \Re (\varphi_n(\omega)))$ 
where $\langle \Re (\varphi_n(\omega)),f\rangle = \Re\langle \varphi_n(\omega),
f\rangle$ $\forall\ f\in C(\tee)$. 
Let $C= \bigcap_k C_k$ and $L = \bigcup_k L_k$.

\noindent {\bf Claim:}
$\mu (L)=0$.

To see the claim, assume that $\mu (L)>0$. Since $L= \bigcup_k L_k$,
 there exists $k\in\nat$ so that 
$\mu (L_k)>0$. 
Consider $\varphi_n^k\in \conv \{\varphi_n,\varphi_{n+1},\ldots\}$ and let 
$\P = \{k\in\nat,\  b(k)>0\}$ and $\N = \{k\in\nat,\  a(k)<0\}$. 
Clearly $\nat = \P\cup \N$.

 Let us assume first that $k\in \P$. 
Using (iii) with $\Gamma = (1,1,1,\ldots)$, for each $\omega\in L_k$, 
there exists $T\in \L(c_0,A)$, $\|T\|\le1$ so that $\Re \langle\varphi_n^{(k)}
(\omega),Te_n\rangle \ge b(k)$. 
Using similar argument as in \cite{RAN4} Lemma 4., one can construct a map 
$T:\Omega\to \L(c_0,A)$ with: 
\begin{enumerate}
\item[a)] $\omega\mapsto T(\omega)e$ is measurable for every $e\in c_0$; 
\item[b)] $\|T(\omega)\|\le 1$ $\forall\ \omega\in \Omega$ and 
$T(\omega)=0$ for $\omega\in\Omega\setminus L_k$. 
\item[c)] $\Re \langle \varphi_n^{(k)} (\omega),T(\omega)e_n\rangle \ge 
b(k)$ $\forall\ \omega\in L_k$.
\end{enumerate}
So we get that 
$$\liminf_{n\to\infty} \int_{L_k} \Re \langle \varphi_n^{(k)}(\omega), 
T(\omega) e_n\rangle\, d\mu(\omega) \ge b(k)\mu(L_k)$$ 
which implies that 
$$\liminf_{n\to\infty} \Big|\int_{L_k}\langle \varphi_n^{(k)}(\omega), 
T(\omega)e_n\rangle\,d\mu(\omega)\Big| \ge b(k)\mu (L_k)\ .$$ 
If $k\in \N$, we repeat the same argument with $\Gamma=(0,0,0,\ldots)$ 
to get that 
$$\liminf_{n\to\infty} \Big|\int_{L_k} \langle \varphi_n^{(k)}(\omega), 
T(\omega)e_n\rangle\, d\mu(\omega)\Big| \ge |a(k)| \mu(L_k)\ .$$ 
So in both cases, if $\delta = \max(b(k)\mu(L_k),|a(k)|\mu(L_k))$, 
there exists a map $T:\Omega\to \L(c_0,A)$ (measurable for the strong 
operator topology) so that 
\begin{equation}\label{star}
\liminf_{n\to\infty} \Big| \int_{L_k}\langle \varphi_n^k (\omega),T(\omega)
e_n\rangle \,d\mu(\omega)\Big| \ge \delta\ .
\end{equation}

To get the contradiction, let 
$$\varphi_n^{(k)} = \sum_{i=p_n}^{q_n} \lambda_i^n |\tilde\sigma (f_i)
(\omega)| = \sum_{i=p_n}^{q_n} \lambda_i^n \tilde\sigma (f_i)(\omega) 
. h_i(\omega)$$ 
with $\sum\limits_{i=p_n}^{q_n} \lambda_i^n =1$, $p_1<q_1<p_2<q_2<\cdots$ 
and $h_i\in L^\infty (\mu,H_\sigma^\infty)$.

Condition \eqref{star} is equivalent to:
$$\liminf_{n \to \infty}\Big| \sum_{i=p_n}^{q_n} \lambda_i^n \int_{L_k}
\langle \tilde\sigma(f_i)(\omega).h_i(\omega), T(\omega)e_n \rangle
\ d\mu(\omega) \Big| \ge \delta.$$
Therefore  there exists $N \in \nat$ so that for each $n \geq N$,
$$\sum_{i=p_n}^{q_n} \lambda_i^n \Big|\int_{L_k}\langle
 \tilde\sigma(f_i)(\omega).h_i(\omega), T(\omega)e_n\rangle \ 
d\mu(\omega) \Big| \geq \delta/2;$$
for each $n \ge N$, choose $i(n) \in [p_n,q_n]$ so that
 $$\Big| \int_{L_k} \langle \tilde\sigma(f_{i(n)})(\omega). h_{i(n)}(\omega),
 T(\omega)e_n \rangle \ d\mu(\omega) \Big| \ge \delta/2 $$ 
and we obtain that for each $n \ge N$,
 \begin{equation}\label{star2}
 \Big|\int_{L_k} \langle \sigma(f_{i(n)}(\omega)), T(\omega)e_n . 
h_{i(n)}(\omega) \rangle \ d\mu(\omega)\Big| \ge \delta/2.
\end{equation}
Notice that for every $\omega \in\Omega$, $T(\omega)e_n \in A$ and 
$h_{i(n)}(\omega) \in H^\infty(\tee)$ so
the product  $T(\omega)e_n . h_{i(n)}(\omega)\in
 H^\infty(\tee)$ and therefore
$$\langle\sigma(f_{i(n)}(\omega)), T(\omega)e_n . h_{i(n)}(\omega)\rangle =
\langle f_{i(n)}(\omega), T(\omega)e_n . h_{i(n)}(\omega)\rangle.$$
For  $n \ge N$, fix 
$$\phi_n(\omega) = \begin{cases}
         T(\omega)e_n . h_{i(n)}(\omega)&\text{$\omega \in L_k$}\\
        0&\text{$\omega \notin L_k$}\end{cases}.$$
If we set $\phi_n =0$ for $n <N$ then
the series $\sum\limits_{i=1}^\infty \phi_i$ is a W.U.C.\ series in $L^\infty (\mu,
H_\sigma^\infty)$: to see this notice that for each $\omega\in\Omega$, 
$\sum\limits_{n=1}^\infty T(\omega)e_n$ is a W.U.C.\ series in $A$ 
(hence in $C(\tee)$) so 
$\sum\limits_{n=1}^\infty |T(\omega)e_n|$ is a W.U.C.\ series in $C(\tee)$. 
Now let $x\in L^1(\mu,\LH)$ (the predual of $L^\infty (\mu,H_\sigma^\infty)$) 
and fix $v\in L^1(\mu,L^1)$ with $\tilde q(v)=x$, 
we have 
\begin{align*}
\sum_{n=1}^\infty |\langle\phi_n,x\rangle| 
&= \sum_{n=1}^\infty |\langle \phi_n,v\rangle| \\
&= \sum_{n=N}^\infty\  
|\langle T(\cdot)e_n. h_{i(n)}(\cdot). \chix_{L_k}(\cdot),v\rangle |\\
&\le  \sum_{n=N}^\infty\  \|h_{i(n)}\| 
\langle |T(\cdot)e_n|,|v|\rangle\\
&\le \sum_{n=1}^\infty \langle |T(\cdot)e_n|,|v|\rangle < \infty\ .
\end{align*}
Now \eqref{star2} is equivalent to: for each $n \ge N$, 
$$ \Big|  \langle \phi_n , 
f_{i(n)}\rangle \Big| \ge \delta/2$$ 
which is a contradiction since $\{f_i,\ i\in \nat\}\subseteq K$ is relatively 
weakly compact and $\sum\limits_{n=1}^\infty \phi_n$ is a W.U.C.\ series. 
The claim is proved. 

To complete the proof of the lemma, let us fix a sequence $(\xi_n)_n$ so that 
$\xi_n\in \conv \{\varphi_n^{(k)},\varphi_{n+1}^{(k)},\ldots,\}$ for every 
$k\in \nat$, we get by (ii) that $\lim\limits_{n\to\infty} \Re \langle \xi_n(\omega),
Te_n\rangle$ exists for every $T\in \L(c_0,A)$; 
we repeat the same argument as above for the imaginary part (starting from 
$(\xi_n)_n$) to get a sequence $(\psi_n)_n$ with $\psi_n\in\conv \{\xi_n,
\xi_{n+1},\ldots\}$ so that $\lim\limits_{n\to\infty} \Im \langle\psi_n(\omega), 
Te_n\rangle$ exists for every $T\in \L(c_0,A)$. 
The lemma is proved.

To finish the proof of the theorem, we will show that for a.e. $\omega$, 
the sequence $(\psi_n(\omega))_{n\ge1}$ is uniformly integrable. 
If not, there would be a measurable subset $\Omega'$ of $\Omega$ with 
$\mu(\Omega')>0$ and $(\psi_n(\omega))_{n\ge1}$ is not uniformly integrable 
for each $\omega\in\Omega'$. 
Hence by Lemma~2, for each $\omega\in\Omega'$, there exists $T\in \L(c_0,A)$ 
so that: 
$$\limsup_{m\to\infty}\ \sup_n |\langle\psi_n(\omega),Te_m\rangle|>0.$$ 
So there would be increasing sequences $(n_j)$ and $(m_j)$ of integers, 
$\delta>0$ so that $|\langle\psi_{n_j}(\omega),Te_{m_j}\rangle|>\delta$ 
$\forall\ j\in \nat$; choose an operator $S:c_0\to c_0$ so that $Se_{n_j} = 
e_{m_j}$; we have $|\langle \psi_{n_j}(\omega),TSe_{n_j}\rangle|>\delta$. 
But by Lemma~3, $\lim\limits_{n\to\infty} |\langle \psi_n(\omega),TSe_n\rangle|$ 
exists so $\lim\limits_{n\to\infty} |\langle \psi_n(\omega),TSe_n\rangle|>\delta$. 
We have  just shown that for each $\omega\in \Omega'$, there exists an 
operator $T\in \L(c_0,A)$ so that $\lim\limits_{n\to\infty} |\langle\psi_n(\omega), 
Te_n\rangle| >0$ and same as before, we can choose the operator $T$ 
measurably, i.e., there exists $T:\Omega\to \L(c_0,A)$, measurable for the 
strong operator topology so that: 
\begin{enumerate}
\item[a)] $\|T(\omega)\|\le 1$ for every $\omega\in\Omega$; 
\item[b)] $\lim_{n\to\infty} |\langle\psi_n(\omega),T(\omega)e_n\rangle|
= \delta(\omega)>0$ for $\omega\in\Omega'$; 
\item[c)] $T(\omega) =0$ for $\omega\notin \Omega'.$
\end{enumerate}
These conditions imply that 
$$\lim_{n\to\infty} \int |\langle \psi_n(\omega),T(\omega)e_n\rangle|\, 
d\mu(\omega) = \int_{\Omega'} \delta(\omega) = \delta >0$$ 
and we can find measurable subsets $(B_n)_n$ so that 
$$\liminf_{n\to\infty} \Big|\int_{B_n}\langle \psi_n(\omega),T(\omega)e_n
\rangle \,d\mu(\omega)\Big| > {\delta\over4}$$ 
and one can get a contradiction using similar construction as in the proof 
of Lemma~3.

We have just shown that for each sequence $(f_n)_n$ in $K$, there  exists a 
sequence $\psi_n \in \text{conv}(|\tilde\sigma(f_n)|,\ |\tilde\sigma(f_{n+1})|,\dots)$
so that for a.e $\omega \in \Omega$, the set 
$\{\psi_n(\omega),\ n\ge 1\}$ is relatively weakly compact in $L^1(\tee)$.
By Ulger's criteria of weak compactness for Bochner space (\cite{U}), the set
$|\tilde\sigma(K)|$ is relatively weakly compact in 
$L^1(\mu,L^1(\tee))= L^1(\Omega \times \tee,\mu \otimes m)$. Hence $\tilde\sigma(K)$ 
is uniformly integrable in $L^1(\Omega \times \tee, \mu \otimes m)$ which is equivalent
to $\tilde\sigma(K)$ is relatively weakly compact in $L^1(\mu,L^1(\tee))$.  
This completes the proof.
\end{pf*} 

Theorem 1. can be extended to the case of spaces of measures. 

\begin{cor}
Let $K$ be a relatively weakly compact subset of $M(\Omega,A^*)$. The
set $\tilde\gamma(K)$ is relatively weakly compact in $M(\Omega,M(\tee))$.
\end{cor}
The following lemma will be used for the proof.
\begin{lem}
Let $\Pi: M(\tee) \to L^1$ be the usual projection. The map $\Pi$ is weak*
to norm universally measurable.
\end{lem}
\begin{pf}
For each $n \in\nat$ and $ 1\leq k <2^n$, let
$D_{n,k}=\{e^{it};\ \frac{(k-1)\pi}{2^{n-1}}\leq t <\frac{k\pi}{2^{n-1}}\}$. 
  Define
for each measure $\lambda$ in $M(\tee)$, 
 $R_n(\lambda)=g_n \in L^1$ be the function 
$\sum\limits_{k=1}^{2^n} 2^n \lambda(D_{n,k})\chi_{D_{n,k}}$.
It is not difficult to see that  the map
$\lambda \mapsto \lambda(D_{n,k})$ is weak*-Borel, so the map $R_n$ is weak* 
Borel measurable as a map from $M(\tee)$ into $L^0$.
 But $R_n(\lambda)$ converges a.e. to the derivative of
$\lambda$ with respect to $m$. If $R(\lambda)$ is such limit, the map $R$
is weak* Borel measurable and therefore $M_s(\tee)=R^{-1}(\{0\})$ is weak*
Borel measurable. Now fix $B$ a Borel measurable subset of $L^1$. Since
$L^1$ is a Polish space and the inclusion map of $L^1$ into $M(\tee)$ is norm
to weak* continious, $B$ is a weak* analytic subset of $M(\tee)$ which implies
that $\Pi^{-1}(B)=B+M_s(\tee)$ is a weak* analytic (and hence weak* 
universally measurable) subset of $M(\tee)$. Thus the proof of the lemma is
complete.
\end{pf}
To prove the corrolary, let $K$ be a relatively weakly compact subset of
$M(\Omega,A^*)$. There exists a measure $\mu$ in $(\Omega,\Sigma)$ so that 
$K$ is uniformly continuous with respect to $\mu$. For each $G \in K$,
choose $\omega \mapsto g(\omega) (\Omega \to A^*)$ a weak*-density of $G$ with 
respect to $\mu$.
Let $g(\omega)=\{g_1(\omega),g_2(\omega)\}$ the unique decomposition of 
$g(\omega)$ in $\LH \oplus_1 M_s(\tee)$. We claim that the function
$\omega \mapsto g_1(\omega)$ belongs to $L^1(\mu,\LH)$. To see this,
notice that 
 the function
 $\omega \to \gamma(g(\omega))=\{\sigma(g_1(\omega)),g_2(\omega)\}$ is a 
 weak*-density of $\tilde\gamma(G)$ with respect to $\mu$. By the above lemma,
 $\omega \mapsto \Pi(\gamma(g(\omega)))=\sigma(g_1(\omega))$
 ($\Omega \to L^1$) is norm measurable
and hence $\omega \mapsto g_1(\omega)$ ($\Omega \to \LH$) is
norm measurable and the claim is proved. 

We get that
 $g(\omega)=\{g_1(\omega) , g_2(\omega)\}$
where $ g_1(.)  \in L^1(\mu,\LH)$ and 
$g_2(.)$ defines a measure in $ M(\Omega, M(\tee))$.
 So $K= K_1 + K_2$ where $K_1$ is a relatively weakly compact subset of 
$L^1(\mu.\LH)$  and $K_2$ is a relatively weakly compact subset of
$M(\Omega,M(\tee))$. It is now easy to check 
$\tilde\gamma(K)= \tilde\sigma(K_1) + K_2$ and an appeal to Theorem 2. 
completes the proof.    
\begin{remark}
Hensgen initiated the study of possible existence and uniqueness of minimum 
norm lifting $\sigma$ from $L^1(X)/H_0^1(X)$ to $L^1(X)$ in \cite{HE}. 
He proved (see Theorem~3.6 of \cite{HE}) that if $X$ is reflexive  
  then 
$\sigma (K)$ is relatively weakly compact in $L^1(X)$ if and only if $K$ is 
relatively weakly compact in $L^1(X)/H_0^1(X)$.  
\end{remark}

\section{The Dunford-Pettis Property}

In this section we prove our main results concerning the spaces $L^1(\mu,\LH)$ 
and $C(\Omega,A)$. 
Let us first recall some characterizations of the Dunford-Pettis property 
that are useful for our purpose. 

\begin{prop}
\cite{D2} 
Each of the following conditions is equivalent to the Dunford-Pettis 
property for a Banach space $X$
\begin{enumerate}
\item[(i)] If $(x_n)_n$ is a weakly Cauchy sequence in $X$ and $(x_n^*)_n$ 
is a weakly null sequence in $X^*$ then $\lim\limits_{n\to\infty} x_n^*(x_n)=0$; 
\item[(ii)] If $(x_n)_n$ is a weakly null sequence in $X$ and 
$(x_n^*)_n$ is a weakly Cauchy sequence in $X^*$ then $\lim\limits_{n\to\infty} 
x_n^*(x_n)=0$.
\end{enumerate}
\end{prop}

It is immediate from the above proposition  that if $X^*$ has the Dunford-Pettis 
property then so does $X$.

We are now ready to present our main theorem. 

\begin{thm} 
Let $\Omega$ be a compact Hausdorff space, the dual of $C(\Omega,A)$ has the
Dunford-Pettis property.
\end{thm}

\begin{pf}
Let $(G_n)_n$ and $(\xi_n)_n$ be weakly null sequences of $M(\Omega,A^*)$ and
$M(\Omega,A^*)^*$ respectively and consider the inclusion map
$J: C(\Omega,A) \to C(\Omega,C(\tee))$. By Corrolary~1, the set
$\{\tilde\gamma(G_n);\ n \in \nat\}$ is relatively weakly compact in 
$M(\Omega,M(\tee))$.

\noindent{\bf Claim}: for each $G \in M(\Omega,A^*)$ and
  $\xi \in M(\Omega,A^*)^*$,
 $\langle G,\xi \rangle =\langle \tilde\gamma(G), J^{**}(\xi)\rangle$.

Notice that the claim is trivially true for $G \in M(\Omega,A^*)$ and 
$f \in C(\Omega,A)$.  For $\xi \in M(\Omega,A^*)^*$, fix a net
$(f_\alpha)_\alpha$ of elements of $C(\Omega,A)$ that converges to 
$\xi$ for the weak*-topology. We have
\begin{align*}
  \langle G,\xi \rangle &= \lim\limits_\alpha \langle G, f_\alpha \rangle \\
  &= \lim\limits_\alpha \langle \tilde\gamma(G), J(f_\alpha) \rangle \\
  &= \langle \tilde\gamma(G), J^{**}(\xi)\rangle
\end{align*}
and the claim is proved.
  
 To complete the proof of the theorem, we use the claim to get that
for each $n \in \nat$, 
  $$\langle G_n, \xi_n \rangle =
\langle \tilde\gamma(G_n),J^{**}(\xi_n)\rangle.$$
Since $(J^{**}(\xi_n))_n $ is a weakly null sequence in $M(\Omega,M(\tee))^*$
and $\{\tilde\gamma(G_n);\ n \in \nat\}$ is relatively weakly compact, we apply
the fact that $M(\Omega,M(\tee))$ has the Dunford-Pettis property (it is an
$L^1$-space) to conclude that
the sequence $(\langle \tilde\gamma(G_n), J^{**}(\xi_n) \rangle)_n$ 
converges to zero and so does the sequence 
$(\langle G_n, \xi_n \rangle)_n$. This completes the proof.
\end{pf} 
  
\begin{cor}
Let $\Omega$ be a compact Hausdorff space and $\mu$ be a finite Borel 
measure on $\Omega$. The following spaces have the Dunford-Pettis prroperty: 
 $L^1(\mu,\LH)$, $L^1(\mu,A^*)$ and $C(\Omega,A)$. 
\end{cor} 

\begin{pf} 
For the space $L^1(\mu,\LH)$, it enough to notice that the space
$L^1(\mu,\LH)$ is complemented in $M(\Omega, \LH)$ which in turn is
 a complemented subspace of $M(\Omega,A^*)$.

For $L^1(\mu,A^*)$, we use the fact that 
 $A^* = \LH \oplus_1 M_S(\tee)$. It is clear that 
$L^1(\mu,A^*) = L^1(\mu,\LH) \oplus_1 L^1(\mu,M_S(\tee))$ and since 
$L^1(\mu,M_S(\tee))$ is an $L^1$-space, the space $L^1(\mu,A^*)$ has the 
Dunford-Pettis property.  
\end{pf} 

\begin{ack}
The author would like to thank  Nigel Kalton for showing him the 
argument used in Lemma~4 and Wolfgang Hensgen for some fruitful comments.
\end{ack}

\bibliography{narciref}
\bibliographystyle{plain}

\end{document}